\documentclass{conm-p-l}

\newtheorem{theorem}{Theorem}[section]

\theoremstyle{definition}

\newtheorem{example}[theorem]{Example}

\theoremstyle{remark}

\numberwithin{equation}{section}

\input psfig.sty

\def\bR{\bf{R}}

\def\bZ{\bf{Z}}

\def\cA{{\mathcal{A}}}

\def\cD{{\mathcal{D}}}

\def\cF{{\mathcal{F}}}
\def\cG{{\mathcal{G}}}
\def\cH{{\mathcal{H}}}

\def\cM{{\mathcal{M}}}

\def\la{\langle}
\def\ra{\rangle}
\def\x{\times}

\def\sk1{\vskip 10pt}

\def\bpm{\begin{pmatrix}}
\def\epm{\end{pmatrix}}
\def\bbm{\begin{bmatrix}}
\def\ebm{\end{bmatrix}}

\def\lb{\label}
\newcommand{\Zn}{{\mathbb Z}}
\newcommand{\fg}{{\mathfrak{g}}}
\newcommand{\fh}{{\mathfrak{h}}}
\newcommand{\cofo}[1]{{\langle #1 \rangle}}
\newcommand{\info}[1]{{( #1 )}}

\begin{document}

\title{Subalgebras of Hyperbolic Kac-Moody Algebras}

\author{Alex J. Feingold}
\address{Dept. of Math. Sci., 
The State University of New York,
Binghamton, New York 13902-6000
}
\email{alex@math.binghamton.edu}

\thanks{AJF wishes to thank the Max-Planck-Institut f\"ur 
Gravitationsphysik, Albert-Einstein-Institut, Potsdam, Germany, for 
support during two stimulating visits.}

\author{Hermann Nicolai}
\address{Max-Planck-Institut f\"ur Gravitationsphysik,
Albert-Einstein-Institut, M\"uhlenberg 1, D-14476 Golm, Germany}
\email{nicolai@aei.mpg.de}

\date{}


\maketitle

\tableofcontents

{\it We dedicate this work to the memory of our friend Peter Slodowy.}

\begin{abstract}
The hyperbolic (and more generally, Lorentzian) Kac-Moody (KM) Lie algebras $\cA$
of rank $r+2 > 2$ are shown to have a rich structure of indefinite KM subalgebras
which can be described by specifying a subset of positive real roots of $\cA$
such that the difference of any two is not a root of $\cA$. Taking these as 
the simple roots of the subalgebra gives a Cartan matrix, generators and 
relations for the subalgebra.  
Applying this to the canonical example of a rank 3 hyperbolic 
KM algebra, $\cF$, we find that $\cF$ contains 
all of the simply laced rank 2 hyperbolics,
as well as an infinite series of indefinite KM subalgebras of rank 3. 
It is shown that $\cA$ also contains
Borcherds algebras, obtained by taking all of the root spaces of
$\cA$ whose roots are in a hyperplane (or any proper subspace). This applies as
well to the case of rank 2 hyperbolics, where the Borcherds algebras have all 
their roots on a line, giving the simplest possible examples.     
\end{abstract}

\section{Introduction} 
More than thirty years after their discovery, indefinite and, more 
specifically, hyperbolic Kac-Moody (KM) algebras remain an unsolved challenge 
of modern mathematics. For instance, there is not a single such algebra
for which the root multiplicities are known in closed form. Even less is
known about the detailed structure of the imaginary root spaces, 
and a manageable explicit realization analogous to the current 
algebra realization of affine algebras appears to be beyond reach 
presently. (However, there are intriguing hints that a ``physical'' 
realization might be found in an extension of Einstein's theory of 
general relativity, see \cite{DHN} and references therein). 
Given this lack of knowledge, 
any new information about these algebras may be of potential value 
and importance.

In this note we point out that hyperbolic KM algebras possess 
a very rich structure of nonisomorphic infinite dimensional subalgebras, 
all of which can be generated by a set of root vectors $E_\alpha$ and
$F_\alpha$, where certain roots $\alpha$ from the root system of 
the original algebra are taken as new simple roots. If these roots
are real, then the subalgebra can be an indefinite KM algebra; however,
we will also exhibit examples where the new simple roots are imaginary,
and the resulting subalgebra is a Borcherds algebra \cite{B1,Jur1,Nie}. Two 
subalgebras are said to be $W$-equivalent if their respective simple roots 
are related by a Weyl transformation from the Weyl group of the original
algebra. Geometrically, some of the non-$W$-equivalent  
subalgebras of indefinite KM algebras can be understood in terms of 
conical sections. Generally there are three ways of slicing 
indefinite (and hyperbolic) KM algebras corresponding to the intersections 
of the light-cone with a hyperplane in the dual of the Cartan subalgebra. 
That dual space is equipped with an indefinite bilinear form, 
$(\alpha,\beta)$, and depending on whether the hyperplane is space-like, 
light-like or time-like (meaning that its normal vector, $\nu$, satisfies 
$\nu^2 = (\nu,\nu) < 0$, $\nu^2 = 0$, or $\nu^2 > 0$) one obtains elliptic, 
parabolic or hyperbolic sections corresponding to a decomposition of the 
original KM algebra with respect to  a finite, an affine or an indefinite 
subalgebra. (The latter need not be hyperbolic, even if the original 
algebra is hyperbolic.) 

Previous attempts to understand indefinite KM algebras were 
almost exclusively based on a decomposition into irreducible 
representations of their affine subalgebras. (Decompositions with respect to 
finite subalgebras seem to have received scant attention so far.) 
Although the theory of representations of affine algebras is reasonably 
well understood, this approach has not yielded much information about 
the behavior of the algebra ``deep within the light-cone''. A complete  
analysis was carried out up to affine level two in \cite{FF} for the 
rank 3 hyperbolic algebra, $\cF$ (also referred to as $HA_1^{(1)}$), 
containing the affine algebra $A_1^{(1)}$. This work was generalized
in \cite{KMW} to the rank 10 hyperbolic KM algebra $E_{10}$, where the
first two levels with respect to the affine subalgebra $E_9 = E_8^{(1)}$ 
were analyzed. This analysis has been extended as far as level 5 
for $\cF$ and also carried out for some levels for some other 
algebras \cite{BKM1,BKM2,BKM3,Ka1,Ka2,Ka3,Ka4,Ka5,KaM1,KaM2}, 
but appears rather hopeless as a method for giving all levels. 
Finally, DDF operators have been used in \cite{GN} to probe some 
root spaces in detail. 

We conclude this introduction with a brief description of the main results 
in each section. In section 2 we give the necessary background and details
about the rank 3 hyperbolic KM algebra, $\cF$, which the most decisive example
for understanding the subject. In section 3 we prove the key theorem which
identifies a class of subalgebras of KM algebras by the choice of a set of
positive real roots whose differences are not roots. We apply that theorem
to the algebra $\cF$ to find inside it all simply laced rank 2 
hyperbolic KM algebras, as well as an infinite series of 
inequivalent rank 3 indefinite KM algebras. We have generalizations of these 
results where the role of the algebra $\cF$ is played by an arbitrary rank $r+2$ 
Lorentzian KM algebra, $\cA$, and we find two infinite series of inequivalent
indefinite KM subalgebras, one series of rank $r+1$ and one series of rank $r+2$. 
In section 4 we go back to the special case of $\cF$, and exploit the 
beautiful geometry of its Weyl group, $W$, which is the hyperbolic triangle group
$T(2,3,\infty) = PGL_2(\bZ)$, through its action on the Poincar\'e disk 
model of the hyperbolic plane. This allows us to find, in principle, a large 
class of KM subalgebras of $\cF$ which are inequivalent to those already found
in section 3. A complete classification of such subalgebras is beyond the 
scope of this paper, but we give a number of representative examples to 
illustrate the idea. We also discuss how the geometry of the Poincar\'e disk 
gives information about the generators and relations defining the Weyl groups
of these subalgebras, and how they are related to the Weyl group $W$. 
We then use this geometrical picture to relate the Weyl groups of the series of
rank 3 indefinite subalgebras of $\cF$ found in section 3 to $W$. 
(Thanks to Tadeusz Januszkiewicz for suggesting and 
explaining this geometrical way of finding subgroups of $W$.) In section 5 
we give a method of finding Borcherds generalized KM algebras inside $\cF$, 
inside and Lorentzian algebra $\cA$, and even inside rank 2 hyperbolics. 

\section{The rank 3 algebra $\cF$}
The canonical example of a hyperbolic KM algebra is the
rank 3 algebra $\cF$ studied in \cite{FF}, whose Cartan matrix is
$$\bmatrix 2 & -1 & 0 \\ -1 & 2 & -2 \\ 0 & -2 & 2 \endbmatrix\ ,$$
whose simple roots are $\alpha_{-1}$, $\alpha_0$, $\alpha_1$, 
and whose Dynkin diagram is
$$\bullet\underset{\hskip -10pt\alpha_{-1}} \ \hbox{-----} \ 
\bullet\underset{\hskip -10pt\alpha_{0}}\   
\vbox{\hbox{-----}\vskip -10pt \hbox{-----}
\vskip -10pt\hbox{-----} \vskip -10pt\hbox{-----}\vskip -3pt}\   
\bullet\underset{\hskip -10pt\alpha_{1}}\ .$$ 
This algebra is the minimal rank hyperbolic KM algebra with both finite and 
affine KM subalgebras. Rank 2 hyperbolic KM algebras are also quite interesting
(\cite{F,KaM2}) because of their connections with real quadratic fields and
Hilbert modular forms \cite{LM}. They have infinitely many $A_1$ subalgebras, 
infinitely many rank 2 hyperbolic KM subalgebras, and our new work shows that
they contain Borcherds subalgebras as well. 
But our main inspiration has been the hyperbolic algebra $\cF$, which seems 
to be the simplest such algebra that incorporates all the essential 
features of higher rank examples. Further evidence that $\cF$ deserves most 
serious study is that it contains all the symmetric rank 2 hyperbolics, as 
well as many other rank 3 hyperbolics. From a physicist's point of view, 
it is attractive because it may appear as a hidden symmetry in (an 
extension of) Einstein's theory of gravity in four space-time 
dimensions (see \cite{DHN} for a review of recent developments
and further references).

Many interesting results about the structure of $\cF$ were  
obtained in \cite{FF}, as well as a relationship with Siegel modular 
forms of genus 2, but a complete understanding
of the root multiplicities has remained elusive despite considerable further
work (\cite{BKM1,Ka3,Ka4,Ka5,KaM1}). 
It is apparent from the Cartan matrix and
Dynkin diagram that $\cF$ contains an affine subalgebra $\cF_0$ of type 
$A_1^{(1)}$ with simple roots $\alpha_0$, $\alpha_1$.
The approach in \cite{FF} is based on the decomposition 
$$\cF = \bigoplus_{n\in\bZ} \cF_n$$ 
with respect to $\cF_0$, that is, with respect to the ``level'' which
is the value of the central element of $\cF_0$. The feature of $\cF$ which 
first attracted the attention of AJF was its Weyl group, $W$, which is a
reflection group 
$$ \la r_{-1},r_0,r_1\ |\ r_{-1}^2 = r_0^2 = r_1^2 = 1, (r_{-1} r_0)^3 = 1, 
(r_{-1} r_1)^2 = 1\ra$$
isomorphic to the hyperbolic triangle group, $T(2,3,\infty)$ and to 
$PGL_2(\bZ)$. The action of $W$ shows that there are infinitely 
many subalgebras of type $A_1^{(1)}$ inside $\cF$, corresponding 
to cusps of the modular group, or to lines of null roots on the 
light-cone. However, these lines of null roots are all related by the Weyl
group, and therefore the corresponding affine subalgebras of $\cF$
are all $W$-equivalent.

To give more details, let $H$ be 
the Cartan subalgebra of $\cF$ and let $H^*$ be its dual space, which has a 
basis consisting of the simple roots $\{\alpha_{-1}, \alpha_0,\alpha_1\}$, and a
symmetric bilinear form whose matrix with respect to this basis is the Cartan
matrix above. That form is indefinite, and can be very conveniently described
as follows. On the space $S_2(\bR)$ of $2 \x 2$ symmetric real matrices define 
a bilinear form by
$$\bpm a & b\\b & c\epm \cdot 
\bpm a' & b'\\b' & c'\epm
= 2bb' - ac' - a'c$$
so the associated quadratic form is 
$$2b^2 - 2ac = -2 \det \bpm a & b\\b & c\epm.$$
Then the root lattice of $\cF$ is isometric to $S_2(\bZ)$ where the entries
$a$, $b$ and $c$ are integers, and the weight lattice of $\cF$ is isometric to
$S'_2(\bZ)$ where $a$ and $c$ are integers but $b\in \frac{1}{2}\bZ$. We
make the correspondence explicit by choosing
$$ \alpha_{-1} = \bpm 1 & 0\\0 & -1\epm, \quad
\alpha_0 = \bpm -1 & -1\\-1 & 0\epm,\quad
\alpha_1 = \bpm 0 & 1\\1 & 0\epm.$$ 
The Weyl group action of $A \in PGL_2(\bZ)$ on 
$N = \bpm a & b\\b & c\epm$ is given by $A(N) = ANA^t$ and the
explicit matrices which represent the three simple reflections are
$$ r_{-1} = \bpm 0 & 1\\1 & 0\epm, \quad
r_0 = \bpm -1 & 1\\0 & 1\epm,\quad
r_1 = \bpm 1 & 0\\0 & -1\epm.$$
The root system of $\cF$ is 
$$\Phi = \{N\in S_2(\bZ)\ |\ \det(N) \geq -1\}$$
which is just the elements $N$ of the root lattice for which the norm squared
is less than or equal to 2, $-2\det(N) \leq 2$. The real roots (Weyl conjugates
of the simple roots) are
$$\Phi_{real} = \{N\in S_2(\bZ)\ |\ \det(N) = -1\}$$
which lie on a single sheeted hyperboloid. 
The light-cone is the set of points $N\in S_2(\bR)$ where $\det(N) = 0$. 
All real roots have multiplicity one, and in the case of $\cF$, the same is true 
of all roots on the light-cone. This comes from $\cF_0$ (with simple roots 
$\alpha_0$ and $\alpha_1$) whose underlying finite dimensional Lie algebra is the 
rank one Lie algebra $sl_2$ of type $A_1$. Any light-cone root, $N\in\Phi$, satisfies 
$\det(N) = 0$ and is $W$ equivalent to a unique one on level 0 of the form 
$\bpm a &0\\0 &0\epm = -a(\alpha_0+\alpha_1)$ with $0\neq a\in\bZ$. The roots $\Phi$
decompose into the disjoint union of positive and negative roots, and we can
distinguish these level  $0$ light-cone roots easily according to whether 
$a$ is negative or positive. 
The roots on level 1 are of the form 
$\bpm a &b\\b &1\epm = -(a+1)(\alpha_0+\alpha_1) + b\alpha_1 - \alpha_{-1}$ 
with $a-b^2\geq -1$ so $a\geq -1$. The affine Weyl group of $\cF_0$ is an 
infinite dihedral group generated by $r_0$ and $r_1$, and it preserves each level. 
Each root on level 1 is equivalent by the affine Weyl group to one of the form 
$\bpm a &0\\0 &1\epm = -(a+1)(\alpha_0+\alpha_1) - \alpha_{-1}$,
and because level 1, $\cF_1$, is a single irreducible $\cF_0$-module (the basic 
module), the multiplicity of that root is a value of the classical partition 
function, $p(a+1)$ \cite{FL}. 
Each root on level 2 is equivalent by the affine Weyl group 
to one of the form $\bpm a &0\\0 &2\epm$ or $\bpm a &1\\1 &2\epm$. 
The $\cF_0$-module structure of level 2 is much more subtle, being 
the antisymmetric tensors in $\cF_1\otimes\cF_1$ with one irreducible 
module removed. If we form a
generating function of the level 2 root multiplicities in those two columns, 
$$\sum_{a=0}^\infty {\rm Mult}\bpm a &1\\1 &2\epm\ q^{2a} + 
{\rm Mult}\bpm a &0\\0 &2\epm\ q^{2a+1} $$
it is remarkable that the first 20 coefficients are again given 
by the partition function $p(\det(N)+1)$. The exact formula 
for that generating function 
$$\left[\sum_{n\geq 0} p(n) q^n\right] \left[\prod_{j\geq 1} (1 - q^{4j-2})\right]
\frac{q^{-3}}{2}\left[ \prod_{j\geq 1} (1 + q^{2j-1}) - \prod_{j\geq 1} (1 - q^{2j-1})
- 2q \right]$$
was one of the main results of \cite{FF}, but further work has shown the 
increasing complexity of higher levels as $\cF_0$-modules. 

There are other ways in which one might naturally decompose $\cF$, for example,
with respect to the finite dimensional subalgebra of type $A_2$ with simple
roots $\alpha_{-1}$ and $\alpha_0$. In that decomposition the analogue of level
of a root $\alpha = \sum a_i \alpha_i$ would be the coefficient $a_1$ of 
$\alpha_1$. Each graded piece of the decomposition would be a finite
dimensional $A_2$ module, and it is not hard to use a computer to find the
irreducible modules which occur and their multiplicities for quite a few 
levels\footnote{The level decomposition with respect to the coefficient 
$a_1$ has been used in \cite{DHN}, and is known up to level $a_1 = 56$ 
(T.~Fischbacher, private communication).}. The question is whether there 
is a useful (recursive) description of those module multiplicities which 
sheds more light on the hyperbolic root multiplicities.

\section{Indefinite subalgebras from subroot systems} 

The main point of the present contribution is that indefinite KM algebras
possess infinitely many non-isomorphic subalgebras, which are of indefinite 
type (in fact, KM subalgebras of equal rank with inequivalent Cartan matrices). 
Unlike finite or affine subalgebras, these subalgebras are themselves not 
yet well understood, but the corresponding decompositions may nevertheless 
provide valuable new viewpoints. 

We will use the following theorem, which allows us to find subalgebras by locating
a set of simple roots for the subalgebra within the root system of the larger
algebra. One cannot just choose an arbitrary subset of positive roots and 
declare them to be the simple roots of a subalgebra. For example, if 
root vectors $E_i$ correspond to simple roots $\beta_i$ and root vectors
$F_j$ correspond to their negatives $-\beta_j$, then one of the required 
Serre relations, $[E_i,F_j] = \delta_{ij} H_i$ could be violated for $i\neq j$ if 
$\beta_i - \beta_j$ were a root of the larger algebra. 

\begin{theorem} \label{subrootsys} Let $\Phi$ be the set of roots of a 
Kac-Moody Lie algebra, $\fg$, with Cartan subalgebra, $\fh$, and let  
$\Phi_{real}^+$ be the positive real roots of $\fg$.   
Let $\beta_1,\cdots,\beta_n\in\Phi_{real}^+$ be chosen such that for all 
$1\leq i \neq j \leq n$, we have $\beta_i - \beta_j\notin \Phi$. For $1\leq i\leq n$
let $0\neq E_i\in\fg_{\beta_i}$ and $0\neq F_i\in\fg_{-\beta_i}$  
be root vectors in the one-dimensional root spaces
corresponding to the positive real roots $\beta_i$, and the negative  
real roots $-\beta_i$, respectively, and let $H_i = [E_i,F_i]\in\fh$. 
Then the Lie subalgebra of $\fg$ generated by $\{E_i,F_i,H_i\ |\ 1\leq i\leq n\}$ 
is a Kac-Moody algebra with Cartan matrix 
$C = [C_{ij}] = [2(\beta_i,\beta_j)/(\beta_j,\beta_j)]$. 
Denote this subalgebra by $\fg(\beta_1,\cdots,\beta_n)$.
\end{theorem}

\begin{proof} By construction, the elements $H_i$ are in the Cartan subalgebra,
$\fh$, so the following relations are clear:
$$[H_i,E_j] = \beta_j(H_i) E_j, \quad
[H_i,F_j] = -\beta_j(H_i) F_j, \quad
[H_i,H_j] = 0.$$ 
Because all of the $\beta_i$ are positive real roots, the matrix $C$ satisfies 
the conditions to be a Cartan matrix. 
For $i\neq j$, the bracket $[E_i,F_j]$ would be in the $\beta_i - \beta_j$ 
root space, but we have chosen the $\beta_i$ such that this difference is not a 
root of $\fg$, so that bracket must be zero, giving the relations
$$[E_i,F_j] = \delta_{ij} H_i.$$
To check that $(ad\ E_i)^{1-C_{ji}}\ E_j = 0$ for all $i\neq j$, it would suffice
to show that $(1 - C_{ji})\beta_i + \beta_j$ is not in $\Phi$. Since $\beta_j$
is a real root of $\fg$, and since $\beta_j - \beta_i\notin\Phi$, the $\beta_i$
root string through $\beta_j$ is of the form
$$\beta_j, \beta_j + \beta_i, \cdots, \beta_j + q\beta_i$$
where $-q = 2(\beta_j,\beta_i)/(\beta_i,\beta_i) = C_{ji}$. Therefore, 
$$\beta_j + (q+1)\beta_i = \beta_j + (1 - C_{ji})\beta_i$$
is not in $\Phi$. The relations $(ad\ F_i)^{1-C_{ji}}\ F_j = 0$ for all $i\neq j$,   
follow immediately since $-(1 - C_{ji})\beta_i - \beta_j$ is not in $\Phi$. This 
shows that all of the Serre relations determined by the Cartan matrix $C$ are 
satisfied by the generators given above. 
\end{proof}

Our first application of this theorem is to 
show that the rank 3 algebra $\cF$ contains all the simply
laced rank 2 hyperbolic KM algebras as subalgebras. The decomposition 
of $\cF$ with respect to its affine algebra $\cF_0$ corresponded to 
slicing with respect to planes parallel to an edge of the light-cone, 
and the affine Weyl group acted on those planes by moving roots 
along parabolas. Decomposition with respect to the finite dimensional 
algebra $A_2$ corresponds to slices which intersect the light-cone 
in circles or ellipses. So it should come as no surprise that there 
should be subalgebras whose decompositions correspond to slices which 
intersect the light-cone in hyperbolas. Consequently we will ``locate''
the rank 2 hyperbolic KM algebras by identifying their simple root 
systems inside the root system of $\cF$. There are two equally good choices
for each algebra, distinguished in the theorem below by a choice of plus
or minus sign. 

\begin{theorem} \label{rank2thm} Fix any integer $m \geq 1$. In the root system
$\Phi$ of $\cF$ we have the positive root vectors
$$\beta_0 = \alpha_{-1} = \bpm 1 & 0\\0 & -1\epm,$$
$$\beta_1 = \beta_1^{(m)} = m(\alpha_0 + \alpha_1) \pm \alpha_1 
= \bpm -m &\pm 1\\\pm 1 &0\epm.$$
Then the KM Lie subalgebra $\fg(\beta_0,\beta_1)$ constructed in 
Theorem \ref{subrootsys} has Cartan matrix 
$$\bbm 2 & -m \\ -m & 2 \ebm$$
which is rank 2 hyperbolic for $m \geq 3$. To indicate the dependence on $m$
we will denote this subalgebra by $\cH(m)$. 
\end{theorem}

\begin{proof} The result follows from Theorem \ref{subrootsys} because 
$\beta_0, \beta_1^{(m)} \in\Phi_{real}^+$ and
$$\beta_1^{(m)} - \beta_0 =  m(\alpha_0 + \alpha_1) \pm \alpha_1 - \alpha_{-1}
= \bpm -m-1 &\pm 1\\\pm 1 &1\epm$$
has determinant $-m-2 < -2$ for $m\geq 1$, so is not a root of $\cF$.  
\end{proof}

While the real root vector $E_0$ of $\cH(m)$ may be taken  
to be the simple root generator $e_{-1}$ of $\cF$, the real  
root vector $E_1$ may be written as the $2m$-fold commutator
$$
E_1 = [e_1,[e_0,[e_1,\cdots  [e_1,[e_0 , e_1] \cdots ]]]
$$
if $\beta_1^{(m)} = m(\alpha_0 + \alpha_1) + \alpha_1$, 
and as the commmutator 
$$
E_1 = [e_0,[e_1,[e_0,\cdots  [e_0,[e_1 , e_0] \cdots ]]]
$$
if $\beta_1^{(m)} = m(\alpha_0 + \alpha_1) - \alpha_1$. 
(Because $\beta_1^{(m)}$ is real, all reorderings of these expressions
are equivalent.)

The above theorem generalizes to any Lorentzian KM algebra, $\cA$, that is,
one which is obtained from an affine algebra by 
the procedure of ``over-extension'', attaching an additional node
by one line only to the affine node in an affine Dynkin diagram. 
For any root $\alpha$ such that $(\alpha,\alpha)\neq 0$, define 
$\alpha^\vee = 2\alpha/(\alpha,\alpha)$. Then the Cartan integers which are
the entries of the Cartan matrix of the rank $r+2$ algebra $\cA$ are given by 
$(\alpha_i,\alpha_j^\vee)$. 

\begin{theorem} \label{rankr1thm}  Let $\{\alpha_i\ | -1\leq i\leq r\}$ be
the simple roots of a Lorentzian KM algebra, $\cA$, so $\alpha_1,\dots,\alpha_r$ 
are the simple roots of a finite dimensional simple Lie algebra, $\alpha_0$
is the affine simple root which generates an affine root system when included
in the list, and $\alpha_{-1}$ satisfies 
$(\alpha_{-1},\alpha_{-1}) = 2$,$(\alpha_{-1},\alpha_0) = -1$, 
$(\alpha_{-1},\alpha_i) = 0$ for $1\leq i\leq r$. Write the affine null root  
$\delta = \sum_{j=0}^r n_j \alpha_j$ where $n_0 =1$. 
Fix any integer $m \geq 0$ and define 
$$\beta_0 = \alpha_{-1}, \quad
\beta_1 = \beta_1^{(m)} = m\delta + \alpha_1, \quad
\beta_j = \alpha_j \;\;\; {\rm for} \; 2\leq j\leq r.$$
Then the KM Lie subalgebra $\cA(\beta_0,\cdots,\beta_r)$ constructed in Theorem 
\ref{subrootsys} has Cartan matrix 
$$C = [(\beta_i,\beta_j^\vee)] = 
\bbm 2 & -m & 0 & \cdots & 0 \\ -m & 2 & * & \cdots & * \\
0 & * & 2 & \cdots & * \\ . &  &  &  &  \\ . &  &  &  &  \\ 0 & * & * & \cdots & 2 
\ebm$$
which is rank $r+1$ of indefinite type for $m\geq 2$, 
where the submatrix obtained by removing the first row and column of $C$ is the 
finite type Cartan matrix $C_r^{fin}$ determined by $\alpha_1,\dots,\alpha_r$.  
We also denote this subalgebra by $\cA_{r+1}^{indef}(m)$. 
\end{theorem}

\begin{proof}
As in Theorem 
\ref{rank2thm}, for $m \geq 0$, the roots $\beta_i$, $0\leq i\leq r$, 
satisfy the conditions of Theorem \ref{subrootsys} and therefore 
determine a KM subalgebra of $\cA$ whose Cartan matrix is as shown because
$(\delta,\delta) = 0$, $(\delta,\alpha_i) = 0$ for $0\leq i\leq r$, and 
$(\delta,\alpha_{-1}) = -1$. Note that the submatrix obtained by removing the
first row and column of $C$ is the finite type Cartan matrix $C_r^{fin}$ 
determined by $\alpha_1,\dots,\alpha_r$. Let $C_{r-1}^{fin}$ be the finite type
Cartan matrix obtained from $C_r^{fin}$ by deleting its first row and column. 
Then we see that
$$\det(C) = 2\ \det(C_r^{fin}) - m^2\ \det(C_{r-1}^{fin}).$$ 
Using the table of values of these determinants on page 53 of \cite{Kac}, we 
find that $\det(C) < 0$ in all cases for $m\geq 2$, 
guaranteeing that the subalgebra will be indefinite.\end{proof}  

For example, if the finite dimensional algebra with simple roots 
$\alpha_1,\dots ,\alpha_r$ is of type $A_r$, then it is easy to see that 
$\det(C) = 2(r+1) - m^2 r$, which will be negative when $m^2 > 2(r+1)/r$. 
For $r\geq 2$, this is true for $m \geq 2$, and for $r = 1$ this is true
for $m\geq 3$. The $r = 1$ case is rather extreme since we only have 
$\beta_0 = \alpha_{-1}$ and $\beta_1 = m\delta + \alpha_1$ giving the rank 2
Cartan matrix already studied in Theorem \ref{rank2thm}.  

The above subalgebras $\cA_{r+1}^{indef}(m)$ of $\cA$ 
for different values of $m\geq 2$ have inequivalent Cartan matrices. 
Therefore, $\cA$ possesses infinitely many nonisomorphic 
indefinite KM subalgebras corresponding to the infinitely many ways
of slicing the forward light-cone in root space by time-like
hyperplanes such that the hyperboloidal intersections
become more and more steeply inclined with increasing $m$.
For rank 3 all subalgebras were again hyperbolic, but the indefinite
KM algebras obtained in this way for higher rank in general are 
no longer hyperbolic even if the original subalgebra is hyperbolic.

The above construction can be further modified as follows to obtain 
subalgebras of rank $r+2$.

\begin{theorem} \label{rankr2thm}  Let the notation be as in Theorem 
\ref{rankr1thm}. Fix any integer $m \geq 1$ and define 
$$\gamma_{-1} = \alpha_{-1}, \quad
\gamma_0 = \gamma_0^{(m)} = (m-1)\delta + \alpha_0, \quad
\gamma_j = \alpha_j \;\;\; {\rm for} \; 1\leq j\leq r.$$
Then the KM Lie subalgebra $\cA(\gamma_{-1},\cdots,\gamma_r)$ constructed in 
Theorem \ref{subrootsys} has Cartan matrix 
$$C = [(\beta_i,\beta_j^\vee)] = 
\bbm 2 & -m & 0 & \cdots & 0 \\ -m & 2 & * & \cdots & * \\
0 & * & 2 & \cdots & * \\ . &  &  &  &  \\ . &  &  &  &  \\ 0 & * & * & \cdots & 2 
\ebm$$
which is rank $r+2$ of indefinite type for $m\geq 1$, 
where the submatrix obtained by removing the first row and column of $C$ is the 
affine type Cartan matrix $C_{r+1}^{aff}$ determined by $\alpha_0,\dots,\alpha_r$.  
We also denote this subalgebra by $\cA_{r+2}^{indef}(m)$. The affine subalgebra
with simple roots $\gamma_0,\cdots,\gamma_r$ has minimal null root 
$\sum_{j=0}^r n_j \gamma_j = m\delta$. 
\end{theorem}

\begin{proof} The proof is as in the previous theorem, but 
the submatrix obtained by removing the first row and column of 
$C$ is the affine type Cartan matrix $C_{r+1}^{aff}$ determined by  
$\alpha_0,\dots,\alpha_r$, and $C_r^{fin}$ is the finite type Cartan matrix  
obtained from $C_{r+1}^{aff}$ by deleting its first row and column, so 
$$\det(C) = 2\ \det(C_{r+1}^{aff}) - m^2\ \det(C_r^{fin}).$$ 
Since $\det(C_{r+1}^{aff}) = 0$, and $\det(C_r^{fin}) > 0$, 
we find that $\det(C) < 0$ in all cases for $m\geq 1$, 
guaranteeing that the subalgebra will be indefinite. Of course, the case
when $m = 1$ just gives the original algebra $\cA$, so the only new content
is for $m\geq 2$. \end{proof}  

For the algebra $\cF$, and for $m \geq 2$, this procedure yields infinitely
many inequivalent rank 3 subalgebras with the Cartan matrices
$$\bbm 2 & -m & 0\\ -m & 2 & -2\\ 0 & -2 & 2 \ebm$$
not previously known to exist inside of $\cF$. In the next section we will
explore the geometrical meaning of these subalgebras using the Weyl group of
$\cF$.

\section{Indefinite subalgebras of $\cF$ from the Weyl group}

A beautiful geometrical way to find subalgebras of $\cF$ is to examine
the action of the Weyl group, $W$, which preserves each surface of 
fixed determinant. 
In fact, it preserves each sheet of the two sheeted hyperboloids with fixed
positive integral determinant. Any one of these sheets can be taken as a model
of the hyperbolic plane, isometric to the Poincar\'e disk, whose boundary 
corresponds to the light-cone. In Figure 1 we have shown the fixed points of
the simple reflections which generate $W$, the triangular fundamental domain 
$\cD$ for $W$, and several other reflection lines corresponding to some other 
positive real roots of $\cF$ which are obtained by applying Weyl group elements to
the simple roots. A complete tesselation of the disk would be obtained by including
all reflection lines, which are in one-to-one correspondence with all of the 
positive real roots of $\cF$, $\Phi^+_{real}$. In that tesselation there would 
be a bijection between the elements of the Weyl group, $W$, and the images of 
the fundamental domain $\cD$. In Figure 1 we see $\cD$ is the triangle with  
angles $\pi/2$, $\pi/3$ and $0$, and sides labelled by the simple roots.  
Six images of $\cD$ form a larger triangle with all vertices on the boundary
and all angles zero. The three sides of that larger triangle are
labelled $\alpha_1$, $r_0 \alpha_1 = \alpha_1 + 2\alpha_0$ and 
$r_{-1} r_0 \alpha_1 = \alpha_1 + 2\alpha_0 + 2\alpha_{-1}$. Since that triangle
has all angles zero, and contains 6 copies of $\cD$, it means that
reflections with respect to these three roots generate a subgroup of $W$ of 
index 6 isomorphic to the hyperbolic triangle group $T(\infty,\infty,\infty)$. 
Using Theorem \ref{subrootsys}, those roots determine a rank 3 hyperbolic subalgebra
of $\cF$, inequivalent to any of the algebras in the series found at the end
of the last section. 

More generally, for any $n\geq 1$, let $S = \{\beta_1,\cdots,\beta_n\}$ be a 
subset of $n$ roots from  $\Phi^+_{real}$, corresponding to 
reflection lines in the Poincar\'e disk, such that $\beta_i - \beta_j\notin\Phi$. 
Then the KM subalgebra
$\cF_S = \cF(\beta_1,\cdots,\beta_n)$ determined by Theorem \ref{subrootsys} has 
Weyl group, $W_S$, generated by the reflections $r_{\beta_1},\cdots,r_{\beta_n}$. 
The Cartan matrix of $\cF_S$ is $C_S = [(\beta_i,\beta_j)]$ since all real roots 
of $\cF$ are of squared length $2$. The fundamental domain $\cD_S$ of $W_S$ is 
a union of images of $\cD$, and the number of such images will be equal to the
index of $W_S$ in $W$. Using Figure 1 we can find choices of $S$ which give 
either finite or infinite index sub-Weyl groups. Our first example was mentioned
above. We leave it to the reader to check the condition that 
$\beta_i - \beta_j\notin\Phi$ for each choice of $S$.  
 
\begin{example} \label{example1} If $S$ consists of the roots 
$$\beta_1 = \alpha_1, \quad \beta_2 = r_0 \alpha_1, \quad 
\beta_3 = r_{-1} r_0 \alpha_1,$$
then $\cF_S$ is a rank 3 subalgebra with Cartan matrix  
$$\bbm 2 & -2 & -2\\-2 & 2 & -2\\-2 & -2 & 2\ebm$$
and Weyl group $W_S = T(\infty,\infty,\infty)$ of index 6 in $W$.
\end{example}

There are other hyperbolic triangles with all angles zero in Figure 1, for
example, the one with sides labelled $\alpha_1$, $\alpha_0$ and 
$r_{-1} r_0 r_1 r_0 \alpha_1 = 3\alpha_1 + 4\alpha_0 + 4\alpha_{-1}$. There are also
6 images of the fundamental triangle in this one, so these three reflections
generate a subgroup of index 6 in $W$. The Cartan matrix generated from these
three roots, taken as simple, is the same as the one above. 
Is there an automorphism
of $\cF$ which interchanges these two isomorphic subalgebras? Are these two
index 6 subgroups of $W$ conjugate? 

\begin{example} \label{example2} If $S$ consists of the roots 
$$\beta_1 = r_1 r_0 \alpha_1 = 3\alpha_1 + 2\alpha_0,$$ 
$$\beta_2 = r_0 r_{-1} r_1 r_0 \alpha_1 = 3\alpha_1 + 6\alpha_0 + 2\alpha_{-1},$$ 
$$\beta_3 = r_{-1} r_0 r_1 r_0 \alpha_1 = 3\alpha_1 + 4\alpha_0 + 4\alpha_{-1},$$ 
then $\cF_S$ is a rank 3 subalgebra with Cartan matrix 
$$\bbm 2 & -10 & -10\\-10 & 2 & -10\\-10 & -10 & 2\ebm$$ 
and Weyl group $W_S = T(\infty,\infty,\infty)$ of infinite index in $W$
because there are no relations between the generating reflections and the area
enclosed by the reflection lines contains an infinite number of copies of $\cD$. 
\end{example}

\begin{example} \label{example3} If $S$ consists of the roots 
$$\beta_1 = \alpha_1,\qquad  
\beta_2 = r_0 r_1 r_{-1} r_0 \alpha_1 = 3\alpha_1 + 6\alpha_0 + 2\alpha_{-1},$$
$$\beta_3 = r_0 r_{-1} r_0 r_1 r_0 \alpha_1 = 3\alpha_1 + 6\alpha_0 + 4\alpha_{-1},$$ 
then $\cF_S$ is a rank 3 subalgebra with Cartan matrix 
$$\bbm 2 & -6 & -6\\-6 & 2 & -2\\-6 & -2 & 2\ebm$$ 
and Weyl group $W_S = T(\infty,\infty,\infty)$ of infinite index in $W$
because there are no relations between the generating reflections and the area
enclosed by the reflection lines contains an infinite number of copies of $\cD$.
\end{example}

\begin{example} \label{example4} If $S$ consists of the roots 
$$\beta_1 = r_0 \alpha_{-1} = \alpha_0 + \alpha_{-1},\qquad 
\beta_2 = r_1 \alpha_0 = 2\alpha_1 + \alpha_0,$$ 
$$\beta_3 = r_0 r_1 r_0 \alpha_1 = 3\alpha_1 + 4\alpha_0,$$ 
then $\cF_S$ is a rank 3 subalgebra with Cartan matrix
$$\bbm 2 & -3 & -2\\-3 & 2 & -2\\-2 & -2 & 2\ebm$$ 
and Weyl group $W_S = T(\infty,\infty,\infty)$ of infinite index in $W$
because there are no relations between the generating reflections and the area
enclosed by the reflection lines contains an infinite number of copies of $\cD$.
\end{example}

\begin{example} \label{example5} If $S$ consists of the four roots 
$$
\begin{array}{lll}
\beta_1 &= r_0 \alpha_{-1} &= \alpha_0 + \alpha_{-1}, \\
\beta_2 &= r_{-1} r_1 r_0 \alpha_1 &= 3\alpha_1 + 2\alpha_0 + 2\alpha_{-1},\\
\beta_3 &= r_1 r_0 \alpha_1 &= 3\alpha_1 + 2\alpha_0, \\
\beta_4 &= r_0 r_1 r_0 \alpha_1 &= 3\alpha_1 + 4\alpha_0
\end{array}
$$
then $\cF_S$ is a rank 3 subalgebra with
Cartan matrix 
$$A = [A_{ij}] = \bbm 2 & -2 & -4 & -2 \\ -2 & 2 & -2 & -10 \\
-4 & -2 & 2 & -2 \\ -2 & -10 & -2 & 2 \ebm$$
and Weyl group 
$$W_S = \la r_{\beta_1},r_{\beta_2},r_{\beta_3},r_{\beta_4}\ |\ r_{\beta_i}^2 = 1\ra.$$ 
>From Figure 1 we see that there are 12 fundamental triangles enclosed by  
these reflecting lines, so the index of $W_S$ in $W$ is 12. \end{example} 

In $\cF_S$, the four $H_i$ are linearly dependent in the three dimensional space 
$H$, as are the new simple roots $\beta_1,\cdots,\beta_4$. In fact, we have
$2\beta_1 - \beta_2 + 2\beta_3 - \beta_4 = 0$. It is certainly possible for a
Cartan matrix to be degenerate and still define a Kac-Moody algebra, but 
Theorem \ref{subrootsys} gives the Serre relations because for 
$1\leq i\neq j\leq 4$, $\beta_i - \beta_j$ is not a root of $\cF$. 
There are several ways to solve the dependency problem (\cite{Kac,Jur1,Nie}),
for example, by adjoining some derivations to the Cartan subalgebra
so as to make the simple roots linearly independent, but then the resulting 
algebra will have Cartan subalgebra larger than the Cartan subalgebra of $\cF$.  
The same considerations occur when generalizing these ideas to the higher rank
algebras $\cA$.

It would be interesting if one could use the geometry of the tesselated disk to 
classify the subgroups of the Weyl group of $\cF$ and use that to classify the
subalgebras coming from Theorem \ref{subrootsys}. 

We would like to finish this section by showing how the series of subalgebras
of $\cF$ found at the end of the last section fit into the geometrical point
of view given in this section. When Theorem \ref{rankr2thm} is applied to 
$\cF$, for $m\geq 1$, the set $S$ is 
$$\gamma_{-1} = \alpha_{-1}, \quad
\gamma_0 = \gamma_0^{(m)} = (m-1)\delta + \alpha_0, \quad
\gamma_1 = \alpha_1$$
so the Weyl group $W_S$ is generated by two of the simple generators of $W$,
$r_{-1}$ and $r_1$, along with one other reflection, $r_{\gamma_0^{(m)}}$.
It is not hard to check that
$$\gamma_0^{(1)} = \alpha_0,\ \gamma_0^{(2)} = r_0 \alpha_1,\  
\gamma_0^{(3)} = r_0 r_1 \alpha_0,\ 
\gamma_0^{(4)} = r_0 r_1 r_0 \alpha_1,\ \cdots,$$ 
which can be seen in Figure 1. The case of $m = 1$ just gives $W$, but when $m = 2$
we see that the three reflecting lines enclose a fundamental domain $\cD_S$
containing three images of $\cD$, and the angles of the triangle are 
$\pi/2$, $0$ and $0$, so in that case $W_S = T(2,\infty,\infty)$ is of index
3 in $W$. But for $m\geq 3$, the three sides do not form a triangle, and there
are infinitely many copies of $\cD$ in their fundamental domain, and the index
of $W_S$ in $W$ is infinite. We would also like to mention that 
$$\bbm -1 & m \\ 0 & 1 \ebm$$
is the matrix in $PGL_2(\bZ)$ which represents the reflection $r_{\gamma_0^{(m)}}$ 
as discussed in section 2.

\section{Borcherds algebras as subalgebras of hyperbolic KM algebras}

Although one might have been surprised to find so many inequivalent indefinite
KM subalgebras in any Lorentzian algebra $\cA$, it is perhaps even more
surprising to find Borcherds algebras as 
proper subalgebras in $\cF$, in any Lorentzian algebra $\cA$, 
and even in rank 2 hyperbolics. As is well known, Borcherds (or generalized 
KM) algebras can be defined in terms of a generalized Cartan 
matrix and a set of generators and relations just like standard
KM algebras, but are distinguished by the existence of
{\em imaginary simple roots}. These correspond to zero or negative
diagonal entries in the Cartan matrix, with corresponding modifications
of the Serre relations \cite{B1}. (See also \cite{Jur1,Nie}.) 
Moreover, the multiplicity of a simple imaginary root may be greater than one.

To explain the basic idea let us return to Theorem \ref{rank2thm}. As we have 
seen, for each $m$, the hyperbolic algebra $\cH(m)$ 
can be embedded into $\cF$ by identifying its two simple roots 
$\beta_0$ and $\beta_1$ in the root system of $\cF$, 
and the corresponding simple root generators as multiple commutators 
of the generators of $\cF$. Let us also choose the $\pm$ sign to be $-$ in
the definition of $\beta_1$. The root system of $\cH(m)$ is contained in 
the linear subspace spanned by $\beta_0$ and $\beta_1$ inside the 
root lattice of $\cF$. As we already explained, this subspace 
gives a hyperbolic section of the light-cone in the root 
lattice of $\cF$. The root space $\cH(m)_\beta$ associated with 
any root $\beta$ in this hyperplane is contained in, but in general 
will not equal the root space $\cF_\beta$. Rather we will have 
strict inequality for ``sufficiently imaginary'' roots $\beta$, viz.
$$
{\rm Mult}_{\cH(m)} (\beta) < {\rm Mult}_{\cF} (\beta).
$$
(Actually the difference in dimension of these two root spaces
will grow exponentially as $\beta$ is moved deeper into the light-cone.) 
As an example, let us take the rank 2 ``Fibonacci'' algebra $\cH(3)$   
(so called because of its connection with the Fibonacci numbers \cite{F}),  
and in it the root $\beta = 2 \beta_0 + 2 \beta_1$.
Then, from the table on page 214 of \cite{Kac}, where this root is denoted 
by $(2,2)$, we have 
$$
{\rm Mult}_{\cH(3)} (\beta) = 1.
$$
On the other hand, with the identification of Theorem \ref{rank2thm},
we have $\beta = 4 \alpha_1 + 6 \alpha_0 + 2 \alpha_{-1}$ as a root
of $\cF$. The Weyl group reflection $r_0$ sends $\beta$ to the root
$4 \alpha_1 + 4 \alpha_0 + 2 \alpha_{-1}$ which is denoted by $(4,4,2)$
in \cite{Kac}, page 215, where the multiplicity of the root is given as $7$, so 
$$
{\rm Mult}_{\cF} (\beta) = 7, 
$$
showing that the root space $\cF_\beta$ contains six independent vectors not 
contained in $H(3)_\beta$.
 
We are thus led to define a new algebra inside $\cF$, as follows. 
Let $\Phi(\cH(m))$ be the set of roots of $\cF$ which are in the plane spanned
by $\beta_0$ and $\beta_1$, so they are the same as the roots of $\cH(m)$. 
Let $\fh_m$ be the span of $H_0$ and $H_1$ from Theorems 
\ref{subrootsys} and \ref{rank2thm}, that 
is, the Cartan subalgebra of $\cH(m)$. Then define the subspace
$$\cG(m) = \fh_m \oplus \bigoplus_{\beta\in\Phi(\cH(m))} \cF_\beta$$  
which is a proper Lie subalgebra of $\cF$ which contains $\cH(m)$ properly.  
The only subtle point involved in checking the closure of $\cG(m)$ under brackets 
is that $[\cF_\beta,\cF_{-\beta}]\subseteq \fh_m$,
but this follows immediately from the invariance of the bilinear form on $\cF$.

We can think of $\cG(m)$ as an extension of $\cH(m)$ in the following way.
For $\beta = a_0\beta_0 + a_1\beta_1 \in\Phi(\cH(m))^+$, $0\leq a_0,a_1\in\bZ$, 
define the $\cH(m)$-height $\cH t(\beta) = a_0 + a_1$. Note that this is not
the same as the height of $\beta$ as a root of $\cF$. Define a sequence of 
extensions, $\cG^{(i)}(m)$, recursively, beginning with 
$\cG^{(1)}(m) = \cH(m)$. For $i>1$ let $\cG^{(i)}(m)$ be the Lie algebra 
containing $\cG^{(i-1)}(m)$ and the additional generators taken from the 
complement of $\cG^{(i-1)}(m)_\beta$ in $\cF_\beta$ for all 
$\beta\in\Phi(\cH(m))^+$ with $\cH t(\beta) = i$. This amounts to adding a 
new imaginary simple root for each such $\beta$ with simple multiplicity
\cite{BGGN,BGN} (not to be confused with the root multiplicity of $\beta$) 
$$\mu(\beta) = \dim(\cF_\beta) - \dim(\cG^{(i-1)}(m)_\beta).$$ 
Then, we have 
$$\cG(m) = \bigcup_{i=1}^\infty \cG^{(i)}(m) .$$
However, to work out the complete system 
of imaginary simple roots, their multiplicities, and the 
associated generalized Cartan matrix even for these simple examples
would be an extremely difficult task. This would in particular 
require full knowledge of the root multiplicities of $\cF$. 
Fortunately, this is not necessary because we can invoke the following 
theorem of Borcherds \cite{B2} (also see \cite{Jur1}).

\begin{theorem} \lb{thm-B1}
A Lie algebra $\fg$ is a Borcherds algebra if it has an almost positive
definite contravariant form $\cofo{\ |\ }$, which means that $\fg$ has
the following properties:
\begin{enumerate}
\item(Grading)\quad
  $\fg=\bigoplus_{n\in\Zn}\fg_n$ with $\dim\fg_n<\infty$ for $n\ne0$;
\item(Involution)\quad
  $\fg$ has an involution $\theta$ which acts as $-1$ on $\fg_0$ and maps
  $\fg_n$ to $\fg_{-n}$;
\item(Invariance)\quad
  $\fg$ carries a symmetric invariant bilinear form $\info{\ |\ }$
  preserved by $\theta$ and such that $\info{\fg_m | \fg_n}=0$ unless
  $m+n=0$;
\item(Positivity)\quad
  The contravariant form $\cofo{x|y}:=-\info{\theta(x)|y}$ is positive
  definite on $\fg_n$ if $n\neq0$.
\end{enumerate}
\end{theorem}

\noindent
Thus, we have

\begin{theorem}\lb{thm-B2}
For all $m\geq 3$, $\cG(m)$ is a Borcherds algebra such that
$$\cH(m) \subset \cG(m) \subset \cF.$$
\end{theorem}

\begin{proof}
All properties listed in Theorem \ref{thm-B1} are satisfied for $\cG(m)$ 
because they are manifestly true for $\cF$. \end{proof}  

In Theorem \ref{rankr2thm} we have seen how any Lorentzian KM algebra 
$\cA$ of rank $r+2$ contains an infinite series of inequivalent indefinite
KM subalgebras, $\cA_{r+1}^{indef}(m)$, for $m\geq 2$. 
Let $\Phi(\cA_{r+1}^{indef}(m))$ be the set of roots of $\cA$ which are in
the hyperplane spanned by the simple roots $\beta_i$, $0\leq i\leq r$, of
$\cA_{r+1}^{indef}(m)$, and let $\fh_m$ be its Cartan subalgebra. Then defining
$$\cG(\cA_{r+1}^{indef}(m)) = 
\fh_m \oplus \bigoplus_{\beta\in\Phi(\cA_{r+1}^{indef}(m))} \cA_\beta$$  
gives a proper Lie subalgebra of $\cA$ generalizing the previous construction.

\begin{theorem}\lb{thm-B3}
For all $m\geq 2$, $\cG(\cA_{r+1}^{indef}(m))$ is a Borcherds algebra such that
$$\cA_{r+1}^{indef}(m) \subset \cG(\cA_{r+1}^{indef}(m)) \subset \cA.$$
\end{theorem}

Let's denote more briefly by $\cG$ any of the 
``hyperplane Borcherds subalgebras" just constructed inside $\cA$, 
and let $\cH$ denote
the indefinite KM subalgebra properly contained in $\cG$. Then we have the
following decomposition 
$$
\cG = \cM_- \oplus \cH \oplus \cM_+  
$$
where $\cM_+$ and $\cM_-$ are supported on positive and negative
roots, respectively. This decomposition corresponds to a similar decomposition
found in \cite{Jur2} for all Borcherds algebras, 
$$
\tilde\cG = \cM_- \oplus (\cH \oplus \tilde\fh) \oplus \cM_+ , 
$$
where $\tilde\fh$ is an infinite dimensional extension of the Cartan subalgebra
of the KM algebra $\cH$ which makes all the 
imaginary simple roots linearly independent.
This extension is analogous to the extension of the Cartan subalgebra mentioned
at the end of the previous section. But this extension would not be contained in
$\cA$. As shown in \cite{Jur2}, $\cM_+$ and $\cM_-$ are free Lie algebras. It 
would be interesting to determine their structure as $\cH$-modules. Similar 
structures were studied in \cite{BGGN}. 

Finally, we would like to note that there are Borcherds algebras inside the
rank 2 hyperbolics, and in particular, in $\cH(3)$, whose positive roots are 
shown in Figure 2. Note that the positive real roots are shown by open circles
and the positive imaginary roots by solid dots. The figure also shows the simple 
reflection lines and root multiplicities. We draw the reader's attention to the 
central vertical line in the figure, and define the subalgebra $\cG$ which is the
direct sum of all the root spaces along that line, including the negative root
spaces not shown and the one-dimensional subspace of the Cartan subalgebra spanned 
by $h_1 + h_2$. This is the simplest example of a Borcherds algebra embedded inside
a hyperbolic KM algebra. In this case we have the decomposition
$$
\cG = \cM_- \oplus sl_2 \oplus \cM_+ 
$$
and it would not be hard to determine the number of free generators in the root
spaces of $\cM_+$ by using the formulas in \cite{Jur2} for the dimensions of 
graded subspaces in free Lie algebras with graded generators.


%
%

\newpage

\section{Appendices}

\centerline{Figure 1: Poincar\'e Disk Model of Hyperbolic Plane Tesselated}
\centerline{By the Hyperbolic Triangle Group $T(2,3,\infty)$}
\vskip 10pt
\centerline{\psfig{figure=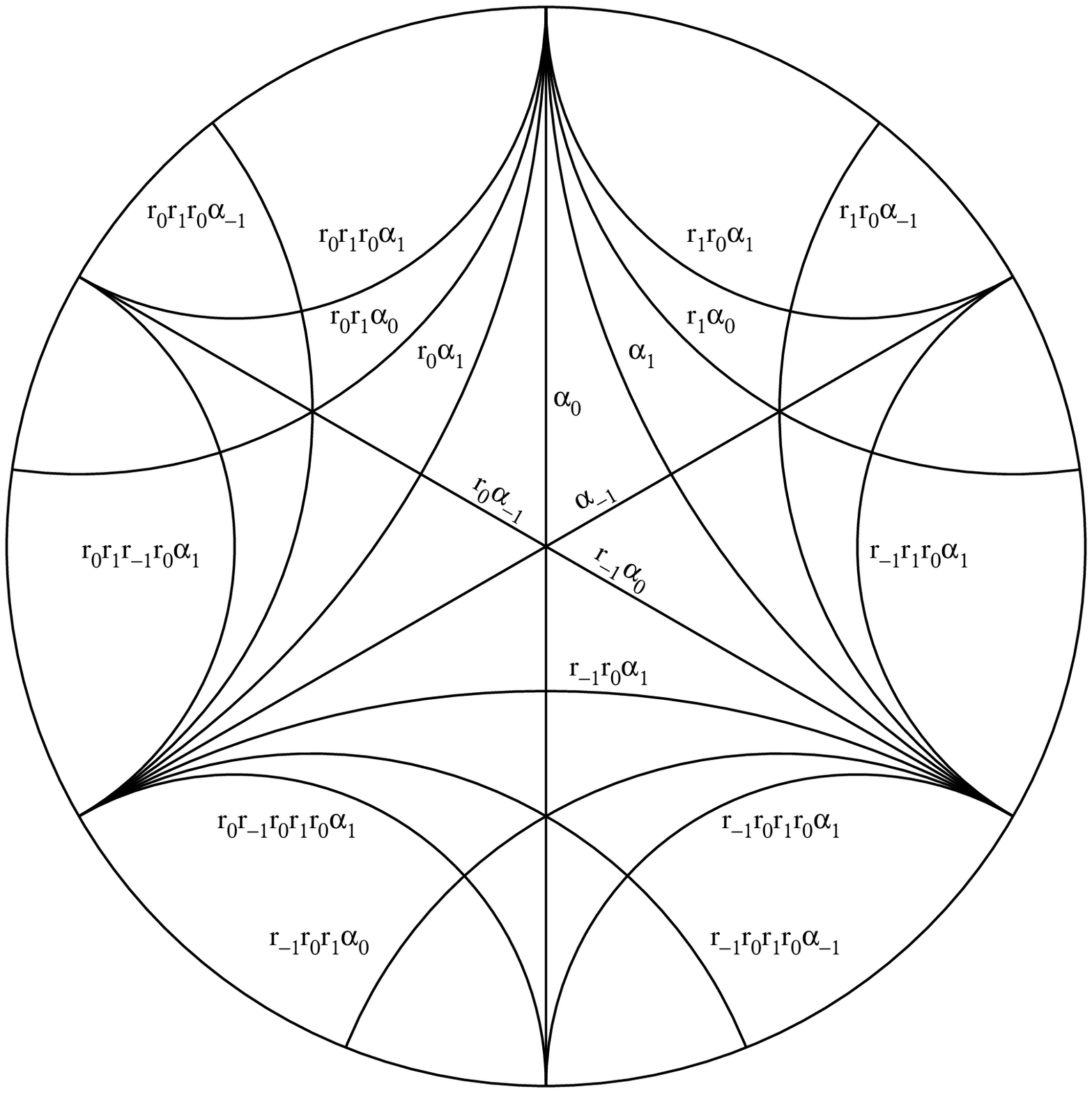,height=6.5in,width=6.5in}}

\newpage
\centerline{Figure 2: Hyperbolic Root System For The Fibonacci Algebra $\cH(3)$}
\vskip 60pt
\centerline{\psfig{figure=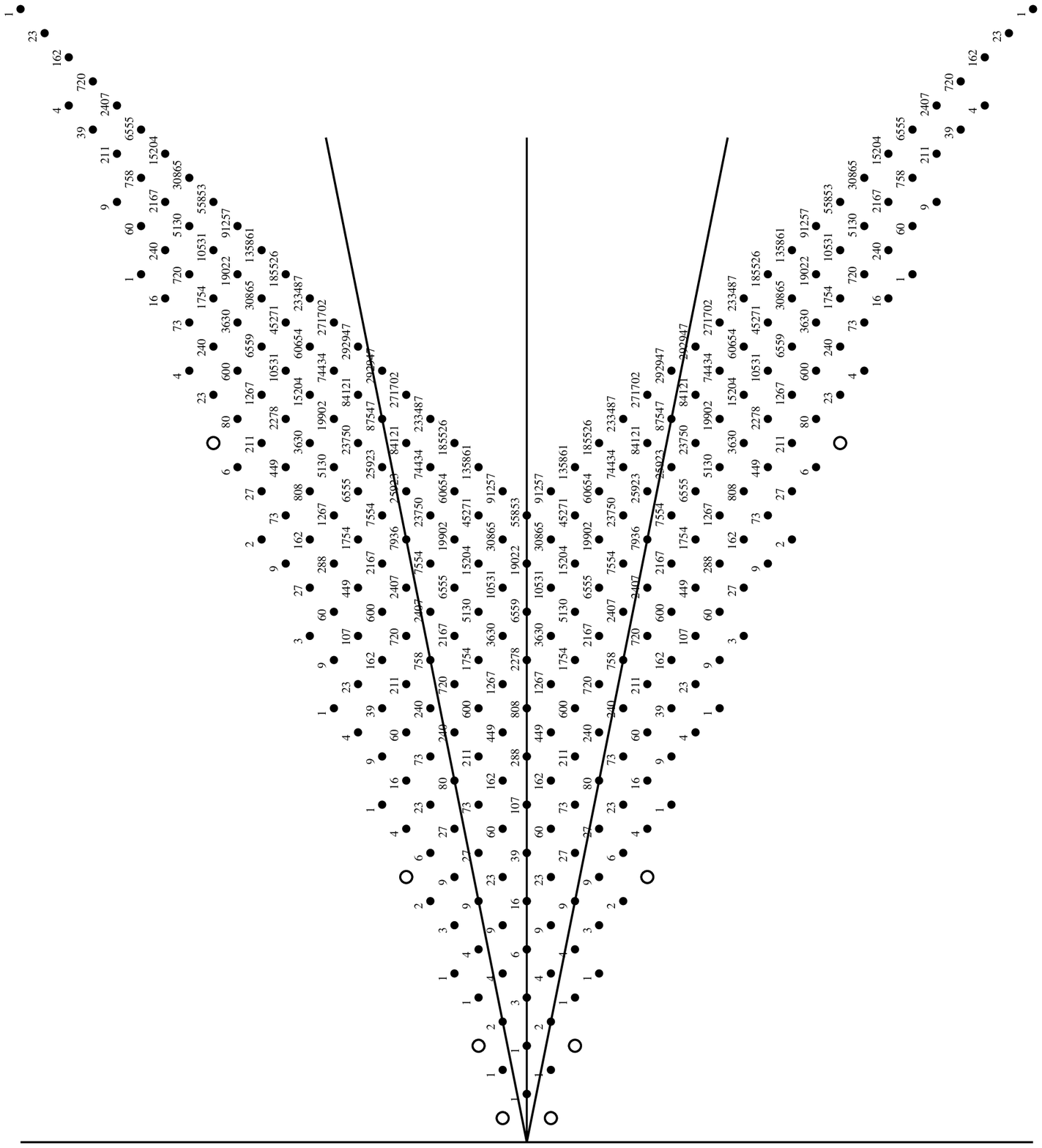,height=7in,width=7in}}

\end{document}